La *Dissertatio de Arte Combinatoria* de Leibniz en seconde lecture

*Alfonso Iommi Echeverría*
*Godofredo Iommi Amunátegui*
Instituto de Física
Universidad Católica de Valparaíso
Casilla 4059, Valparaíso-Chile

Abstract

This article considers the *Dissertatio de Arte Combinatoria*, published in 1666 and relatively neglected by Leibniz's scholars. However in recent times the tide seems to be changing. Our work presents three main parts (sections II, III and IV, respectively). The first one exposes the contributions of M. Serres and E. Knobloch, viewed from our standpoint. The second one emphazises the role played by the symmetric group as the underlying mathematical structure of the Opus. Finally we treat the *caput variationis* which, in our opinion, is the fundamental concept of the Dissertatio. Moreover, it is suggested that the importance of such a theoretical insight may be traced in many aspects of Leibniz's subsequent thought.

# I

La *Dissertatio*[1] date de 1666. Par la suite, Leibniz fait souvent allusion à cette œuvre de sa prime jeunesse. Avec un regard critique, mais sans pour autant méconnaître cette éclosion initiale de sa pensée. Ainsi, à l'occasion d'une réimpression faite sans son autorisation en 1690, il montre du doigt la faiblesse structurelle du livre: *oeconomia operis, in qua multa possent mutari in melius*[2]. Il en parle plus longuement dans les *Nouveaux essais*: "C'est ce qui est arrivé à mon *Art des Combinaisons*, comme je m'en suis déjà plaint. C'était un fruit de ma première adolescence, et cependant on le réimprima longtemps après sans me consulter et sans marquer même que c'était une seconde édition, ce qui fit croire à quelques-uns, à mon préjudice, que j'étais capable de publier une telle pièce dans[sic] un âge avancé ; car quoiqu'il y ait des pensées de quelque conséquence, que j'approuve encore, il y en avait pourtant aussi qui ne pouvaient convenir qu'à un jeune étudiant"[3]. Nous n'avons certes pas l'intention de faire ce partage dans la *Dissertatio* elle-même ni a fortiori dans la pensée ultérieure du philosophe. Tout au plus pourrait-on interroger un texte publié par Couturat[4] où il semble que Leibniz se soit plutôt tourné du côté de la Caractéristique pour y trouver des "conséquences": "Il y a plus de 20 ans... que je m'avisa d'une méthode qui nous mène infailliblement à l'analyse générale des connaissances humaines... comme on peut juger par un petit traité que je fis imprimer alors, où il y a quelques choses qui sentent le jeune homme et l'apprentif [sic] mais le fond est bon... Je trouva donc qu'il y a certains Termes primitifs... lesquels estant [sic] constitués, tous les raisonnements se pourraient déterminer à la façon des nombres". Après Couturat, maints auteurs ont considéré la *Dissertatio*. Pourtant les pages qui lui sont consacrées demeurent rares en regard de la littérature publiée. Nous porterons notre choix sur deux lectures récentes: celle de Serres[5] –philosophique– et celle de Knobloch –mathématique[6].

---

[1] G. W. Leibniz: "Sämtliche Schriften und Briefe", hrsg. v. Deutsche (à présent Berlin-Brandenburgische) Akademie der Wissenschaften (ci-après nous désignerons, suivant l'usage, cette édition par A), Vol VI, 1, pp.163-230. Dans ce travail les passages en français de la *Dissertatio* ont été traduits par nous.

[2] G. W. Leibniz, "Sämtliche Schriften und Briefe", hrsg. v.d. Preuss (à présent Deutschen) Akad. d. Wiss., Darmstadt- Leipzig-Berlin, depuis 1923, ser. VI, vol.2, pp. 549-550.

[3] G. W. Leibniz: "Nouveaux essais sur l'entendement humain", chronologie, notes et introduction par J. Brunschwig, Garnier-Flammarion, Paris, 1966, p. 338.

[4] L. Couturat: "Opuscules et fragments inédits de Leibniz", G. Olms, Hildesheim, 1961, pp. 175-182. Couturat signale que cet opuscule date, au plus tôt, de 1686.

[5] M. Serres: "Le système de Leibniz et ses modèles mathématique", P.U.F., 2ème **éd**. 1982.

[6] E. Knobloch: "Die mathematischen Studien von G. W. Leibniz zur Kombinatorik", Wiesbaden, 1973 (ci-après Studien); "Der Beginn der Determinantentheorie", Hildesheim, 1980 (ci-après Determinantentheorie); "Déterminants et élimination chez Leibniz", Revue d'Histoire des Sciences, 2001, 54/2, pp. 143-164 (ci-après

Le présent article prend ces travaux en compte et reprend la réflexion, mais on ne peut ajouter "là où ils l'ont laissé", car notre interprétation ne coïncide pas toujours avec les leurs.

Après en avoir esquissé les contours, nous procéderons en deux étapes. Tout d'abord nous essayerons de montrer que l'*Ars Combinatoria* peut être envisagée à partir d'une structure algébrique –inconnue de Leibniz–, le groupe des permutations. Pour justifier ce véritable pari théorique, nous devrons introduire les définitions et les rudiments mathématiques nécessaires. Puis nous fixerons notre attention sur la notion de "*caput variationis*". Son rôle est, nous semble-t-il, crucial. Nous verrons qu'il est possible d'en faire le fondement même de l'œuvre de Leibniz.

## II

Une lecture souple et attentive du livre de M. Serres permet d'y percevoir d'emblée l'importance du *caput variationis* (qu'il traduit par "facteur invariant"). Pour en suivre le fil, nous pouvons commencer par cette phrase: "Une autre expression y revient souvent: *caput variationis*, moins importante pour la mathématique, mais, à mon avis capitale pour la philosophie" (p. 413). Nous reviendrons sur la portée mathématique du concept. Penchons-nous tout d'abord sur l'horizon philosophique qui en résulte.

Sans considérer les pages dédiées à la *Dissertatio* (pp. 409–442), on le retrouve de manière implicite ou explicite –le plus souvent d'ailleurs– tout au long de l'ouvrage. Soulignons ces occurrences sur la fin: "il n'est pas d'invariant dans une variation qui ne varie enfin pour un autre invariant, le *caput* devient note pour un nouveau *caput*" (p. 740) et "Dans l'espace combinatoire, nous tournons *ad libitum* autour d'autant de *capita* qu'il y a d'éléments" (p. 791).

Quelles raisons portent Serres à assigner au *caput variationis* un statut privilégié dans son commentaire ? À notre avis la réponse à cette question est à la fois évidente et cachée. Revenons au texte de Serres. Page 102 nous lisons: "Tout se passe comme si chaque point terminal d'une voie méthodique ... pouvait être le point de départ d'une nouvelle structure en étoile autour d'un concept unique de la philosophie et vers la

Déterminants 2001);"The Mathematical studies of G.W. Leibniz on Combinatorics", Historia Mathematica 1, 1974, pp. 409-430 (ci-après: On Combinatorics 1974); "Unbekannte Studien von Leibniz zur Eliminations- und Explikationstheorie", Archive for History of Exact Sciences, vol. 12, n°2, 1974, pp. 142-173 (ci-après Eliminationstheorie 1974); "Der Leibnizsche Dialogus de arte computandi", NTM Schriftenreihe zur Geschichte der Naturwissenschaft, Technik, Medezin, 11, 1974, pp. 10-32; "Studien von Leibniz zum Determinantenkalkül", Studia Leibnitiana Supplementa, 13, 1974, pp. 37-45.

multiplicité des concepts de la philosophie. Autrement dit, chaque concept peut être considéré comme un *caput variationis*, au sens du *De Arte*, point fixe qui ordonne une permutation. Tout concept est une complexion des concepts du système. L'organisation ou l'architectonique de la philosophie leibnizienne n'est pas fondée sur un ordre univoque ou irréversible, mais sur la *variationum* ou *mutatio relationis*. Ainsi l'espace du système, partout centré, répète encore la situation globale de la communication monadique: au pluralisme ontologique des âmes percevantes répond le pluralisme architectonique des centres universels, des parties totales, des concepts comme *capita variatiorum*[22].

Cette conception du *caput* établit une sorte de clé de composition de l'œuvre de Serres tant elle en oriente le sens et en détermine la forme. Par son truchement se noue une correspondance étroite entre le système leibnizien et la structure du discours qui l'expose. Son importance paraît bien assise. Et ce non seulement dans le cadre de cette interprétation de la pensée de Leibniz. Plus loin nous essayerons d'en faire le noyau conceptuel de la *Dissertatio*.

Pourtant nous avons vu que, d'après le même auteur, le "facteur invariant" serait "moins important pour la mathématique". À titre de conclusion, il signale que: "L'apport du *De Arte* à la science Mathématique se résume à cela, qui peut paraître court, mais qui n'a pas été essentiellement amélioré par l'algèbre classique, jusqu'à l'apparition de la théorie des groupes" (p. 415). Par après, nous traiterons cet aspect de la *Dissertatio* et, implicitement nous discuterons cette conclusion. Remarquons, de suite, que le concept de groupe hante, pour ainsi dire, maint passage du livre de Serres. Au moins en deux occasions son apparition est littérale: "L'*Essay de Dynamique* nous avait amené à l'idée que le terme "harmonie" désignait confusément un tableau d'éléments et d'opérations fermé sur soi (la page de garde du *De Arte* est porteuse de la même intuition). Ici l'intuition se précise: la table est harmonique parce qu'elle est structurée comme un groupe" (p. 531). De même: "La vérité n'est pas ici ou là, et n'est pas en tel lieu, exclusivement à tel autre, elle est constituée par l'ensemble, par le tableau, par le groupe des substitutions ... la vérité est l'invariant dans la substitution des identiques ..." (p. 628).

Revenons un instant à cette perspective du *caput variationis* qui en restreint l'horizon mathématique. Ici, il nous semble qu'un point sensible sinon critique est touché. Serres lui–même remarque: "Notons qu'il arrive souvent que ce soit dans les thèses de philosophie qu'on trouve les structures mathématiques les plus fortes" (p. 47). G.-G. Granger[7] abonde dans cette direction: Leibniz présente "un des très rares exemples d'une création mathématique ... novatrice ... associée dès son origine ... à des vues logiques et métaphysiques où elle trouve son impulsion initiale et l'orientation de son mouvement". Partant, le *caput*, à moins de faire exception, devrait avoir une importance mathématique comparable à sa portée philosophique.

---

[7] G.-G. Granger, *Philosophie et mathématique leibniziennes*, in *Revue de Métaphysique et de Morale*, 1981, n°1, pp. 1–37.

Dans les deux sections suivantes de cet article[8], nous traiterons ces aspects de la *Dissertatio*. On y trouvera une discussion, implicite, des vues exposées plus haut.

E. Knobloch soutient que l'art d'inventer fut de loin la préoccupation la plus importante de Leibniz pendant toute sa vie. Cet art s'appuie sur l'art caractéristique (choix de signes) et sur l'art combinatoire. Ce dernier, pris en son sens leibnizien, comprend un domaine mathématique fort étendu.

Knobloch en énumère cinq parties: la combinatoire élémentaire, les fonctions symétriques, les partitions, les déterminants et la théorie des probabilités. Cet auteur a étudié les manuscrits ainsi que les publications portant sur ces thèmes et leur a consacré – notamment aux quatre premiers- divers ouvrages et travaux (voir note 6). Dans "Studien" toute une section du premier chapitre est centrée sur la *Dissertatio* (pp.23-53). L'importance du *caput variationis* y est soulignée à plusieurs reprises (en particulier pp.30-32, pp.47-50). Un article –"On Combinatorics 1974"- condense en quelque sorte ce passage du livre. Notons-y en premier lieu cette remarque: "Les études les plus récentes sur l'*Ars combinatoria* dues à Kutlumuratov (1964) et Michel Serres (1968) ne vont guère au–delà des descriptions antérieures en ce qui concerne les aspects mathématiques" (p. 410). De même, sur le Problème qui traite du *caput variationis*: "Le Problème 7 tient compte des permutations qui contiennent un *caput* c'est-à-dire un sous-ensemble projeté sur lui-même par une permutation ou bien qui demeure invariant (cas spécial)" (p. 415). Plus loin l'approche de l'*Ars Combinatoria* est définie par l'expression "caput–theory". Pour nous, certes, cette idée est centrale.

Knobloch signale l'apport de Leibniz quant aux partitions ("Studien", chap.3; "On Combinatorics 1974", pp.417-426) ainsi que ses contributions à la théorie des permutations cycliques ("Studien", pp.44-45 et "Determinantentheorie", p. 44). Mais nulle part il ne mentionne une relation entre la théorie des groupes et le texte du jeune philosophe.

---

[8] Quelques-unes des thèses ici développées ont été communiquées au VI Internationale**r** Leibniz Kongress (I.L.K.) : G. Iommi Amunátegui, "*Dissertatio de Arte Combinatoria*: A Group-theoretical approach", hsgb. von H. Breger, Vorträge I. Teil pp. 9-16, Hannover 1994 et au VII I.L.K. , G. Iommi Amunátegui et A. Iommi Echeverría: "Le Concept de *Caput Variationis* dans la *Dissertatio* de Leibniz –esquisse d'une théorie", hsgb. von H. Poser, Vorträge 2. Teil pp. 530-535, Berlin, 2001.

**III**

Ci-après nous nous proposons d'établir quelques liens entre la structure de la *Dissertatio* et celle du groupe symétrique $S_n$. Ce propos semble à première vue quelque peu surprenant sinon anachronique, puisque l'oeuvre de Leibniz vit le jour presque deux siècles avant l'introduction formelle de la théorie des Groupes en Mathématique. On ne trouve nulle occurrence explicite du concept de "groupe" –ni a fortiori du terme lui-même– dans le texte. Toutefois l'analyse de divers passages permet d'y déceler une sorte d'action de présence de cette structure algébrique. Il nous paraît même possible d'avancer que ce groupe joue le rôle d'un "principe architectonique" dans le *De Arte*.

<u>Le Contexte mathématique</u>

Soit un ensemble G –fini ou infini– d'éléments a, b, c, ... muni d'une loi de composition qui associe deux de ces éléments à un troisième élément de G. Il existe un élément e de G, l'identité, pourvu de la propriété ke = ek = k ainsi que $k^{-1}$ tel que $kk^{-1} = k^{-1}k = e$; $k^{-1}$ est l'inverse de k, k et $k^{-1}$ appartiennent à G. L'ensemble G est alors un groupe. Un sous-ensemble H de G pourvu de la loi de composition opérant dans G et qui est un groupe est un sous–groupe de G. Un groupe qui a un nombre fini d'éléments est un groupe fini et le nombre de ses éléments est l'ordre du groupe. Si l'on fixe un élément g et que l'on considère un autre élément p qui varie parmi les éléments du groupe, g et $gpg^{-1}$ sont appelés éléments conjugués: ils appartiennent à la même classe de conjugaison.

Le groupe dont les éléments sont les n! permutations de n symboles (1, ..., n) est le groupe symétrique Sn. Les éléments de Sn sont, en fait, les opérations définies par le passage d'un arrangement type, par exemple (1, ..., n) à un arrangement quelconque de ces n symboles. Une permutation de Sn peut se formuler en écrivant sous la disposition standard le nouvel arrangement des symboles. Par exemple:

$$S_3 = \left\{ \begin{pmatrix} 1\,2\,3 \\ 1\,2\,3 \end{pmatrix} \begin{pmatrix} 1\,2\,3 \\ 1\,3\,2 \end{pmatrix} \begin{pmatrix} 1\,2\,3 \\ 2\,1\,3 \end{pmatrix} \begin{pmatrix} 1\,2\,3 \\ 2\,3\,1 \end{pmatrix} \begin{pmatrix} 1\,2\,3 \\ 3\,1\,2 \end{pmatrix} \begin{pmatrix} 1\,2\,3 \\ 3\,2\,1 \end{pmatrix} \right\}$$

Une permutation de la forme $\begin{pmatrix} 1\,2\,...\,r\text{-}1\,r+1\,...\,n \\ 2\,3\,...\,r \quad r+1\,...\,n \end{pmatrix}$ est appelée cyclique ou bien un r–cycle. Une notation plus succincte est (1, 2, ..., r) (r+1) ... (n) où les symboles qui subissent une permutation cyclique sont mis ensemble entre parenthèses. Les permutations qui échangent deux symboles sont appelées transpositions. Notons qu'un 1–cycle signifie que le symbole correspondant conserve sa position initiale. Dans le cas de $S_3$, nous avons donc:

$$S_3 = \{(1)(2)(3); (1)(23); (3)(12); (123); (132); (2)(13)\}$$

Si $\alpha_1$ est le nombre de 1–cycles, $\alpha_2$ le nombre de 2–cycles etc pour toute permutation, il vient:

$$\alpha_1 + 2\alpha_2 + ... + n\,\alpha_n = n$$

Ainsi pour la permutation de $S_6$ $\begin{pmatrix} 1\,2\,3\,4\,5\,6 \\ 1\,4\,3\,6\,5\,2 \end{pmatrix} = (1)(3)(5)(246)$

l'on a:

$$\alpha_1 = 3 \quad \alpha_3 = 1 \quad \alpha_2 = \alpha_4 = \alpha_5 = \alpha_6 = 0$$

donc

$$3 + 2\cdot 0 + 3\cdot 1 + 4\cdot 0 + 5\cdot 0 + 6\cdot 0 = 6.$$

Une permutation qui se décompose en $\alpha_1$ 1–cycles, $\alpha_2$ 2–cycles etc possède une structure cyclique ($\alpha$). On dit que l'ensemble des permutations ayant la même structure ($\alpha$) forme la classe ($\alpha$) des permutations de n symboles. Le nombre des permutations qui appartiennent à la classe ($\alpha$) est donné par la formule de Cauchy:

$$\text{Ordre de la classe } (\alpha) = \frac{n!}{1^{\alpha_1}\alpha_1!\,2^{\alpha_2}\alpha_2!\,...\,n^{\alpha_n}\alpha_n!}$$

(avec 0!=1)

## Caput variationis[*]

Le dernier paragraphe du problème VII nous introduit, si l'on peut dire, dans l'atelier de Leibniz: *Hoc problema casuum multitudo operosissimum effecit, ejusque nobis solutio multo et labore et tempore constitit*. Beaucoup de travail et de temps furent consacrés à le résoudre. Pourtant la phrase suivante nous concerne encore davantage: *Sed aliter sequentia problemata ex artis principiis nemo solvet*. Cette solution donc donne la possibilité de résoudre la série des problèmes posés dans la *Dissertatio* et, partant, de comprendre l'ordre de leur séquence *ex artis principiis*. Les principes de l'Art doivent être cherchés ici même dans la formulation de ce septième problème: "*Dato capite variationes reperire*"[9].

---

[*] Les explications de certains termes ont èté rassemblés dans un Glossaire en fin d'article. De même, il nous faut souligner que nos traductions du latin sont aussi des «interprétation».

Pour saisir le sens de l'expression *caput variationis* examinons-en la définition: *Caput variationis est positio certarum partium*[10]: le facteur invariant d'une variation est la position de ses parties fixes. Il s'agit d'un "invariant" c'est-à-dire d'un élément qu'une variation laisse à sa place. Remarquons qu'ici "variation" est synonyme de "permutation". Or nous savons que les permutations de n symboles constituent le groupe symétrique Sn. Il nous faut tenir compte des parties d'un arrangement de n éléments qui demeurent inchangées sous l'action d'une transformation de Sn.

Ailleurs Leibniz écrit explicitement les 24 permutations des 4 lettres a, b, c, d et il le fait en fixant à tour de rôle chacun de ces symboles. Par exemple, si le "facteur invariant" est a:

(1)  a   b   c   d
(2)  ·   ·   d   c
(3)  ·   c   b   d
(4)  ·   ·   d   b
(5)  ·   d   b   c
(6)  ·   ·   c   b

Les points indiquent la lettre qui occupe la même place à la ligne précédente: *Puncta significant rem praecedentis lineae directe supra positam*[11]. Moyennant la notation cyclique nous pouvons écrire:

(1)  (·)   (·)   (·)   (·)         $\alpha_1 = 4$
(2)  (·)   (·)   (··)              $\alpha_1 = 2$   $\alpha_2 = 1$
(3)  (·)   (·)   (··)              $\alpha_1 = 2$   $\alpha_2 = 1$
(4)  (·)   (···)                   $\alpha_1 = 1$   $\alpha_3 = 1$
(5)  (·)   (···)                   $\alpha_1 = 1$   $\alpha_3 = 1$
(6)  (·)   (·)   (··)              $\alpha_1 = 2$   $\alpha_2 = 1$

---

[9] G.W. Leibniz, A, VI, 1, pp.219-220.
[10] G. W. Leibniz, A, VI, 1, p.173.
[11] G.W. Leibniz, A, VI, 1, p.212.

Un 1–cycle signifie qu'une lettre demeure à sa place. C'est le cas pour le facteur invariant a dans les six permutations écrites ci–dessus. Mais en outre d'autres symboles sont inchangés en certaines d'entre elles: dans la permutation (1), quatre lettres; dans les permutations (2), (3) et (6) deux lettres. Seulement les cas (4) et (5) satisfont à la condition: "la lettre a –exclusivement– est fixée". Notre interprétation, basée sur la théorie des groupes, n'est donc juste que partiellement. Il n'empêche que même avec cette restriction, à notre avis, on peut signaler que:

a) le concept de "facteur invariant" prend un sens précis, issu de la théorie des groupes.
b) Leibniz est tout à fait clairvoyant quant à son importance: ce problème VII montre toute sa richesse théorique dans le calcul combinatoire développé dans la *Dissertatio*.
c) Lors de la définition du concept, l'auteur ajoute: "v. Infra prob. 7". C'est dire qu'une notion–clef est mise en rapport avec un problème qui, à son tour, est établi sur une base algorithmique.

La variation relative de la position

Voici l'énoncé du cinquième problème: "*Dato numero rerum variationem situs mere relati seu vicinitatis invenire*"[12] ("Trouver la variation relative de la position ou du voisinage d'un nombre donné de choses"). Pour en comprendre la portée l'on fera appel à certaines définitions: le lieu (ou place) est la position des parties (*situs est localitas partium*) et il peut être absolu (les parties avec le tout) ou relatif (les parties avec les parties). La variation signifie un changement de relation (*variatio est mutatio relationis*).

Le lieu (place, endroit) absolu peut se représenter par une ligne droite dans laquelle le nombre des lieux et leurs distances par rapport aux deux points extrêmes de la ligne (*initio et fine*) doivent être considérés.

Le lieu (place, endroit) relatif est illustré par un cercle dans lequel ce n'est point l'ordre des lieux qui importe mais leur voisinage. Ainsi l'ordre des dispositions suivantes –

---

[12] G.W. Leibniz, A, VI, 1, p.217.

abcd, bcda, cdab, dabc– diffère mais leur place relative est, pour ainsi dire, en cet unique voisinage:

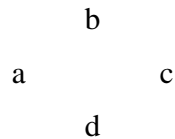

Le nombre des variations de l'ordre de quatre choses est 24 (*Res 4 possunt transponi modis 24*). Le problème admet une nouvelle formulation: "Trouver le nombre de voisinages d'un certain nombre de choses".

Nulle solution générale ne se trouve dans le texte mais un exemple, celui que l'on a exposé plus haut: abcd, bcda, cdab, dabc correspond à une variation écrite circulairement. Donc, les 24 variations doivent être divisées par le nombre de choses –en ce cas, 4– et la solution est:

$$\frac{24}{4} = 6$$

Or 6 est le nombre des variations de trois choses (6=3·2·1=3!). Aussi, pour n choses, le nombre des variations relatives de la position est donné par:

$$\frac{n!}{n} = (n-1)!$$

Il se trouve que chacune des notions considérées peut être traduite –presque mot par mot– en termes de la théorie des groupes. Le cercle qui représente la place relative équivaut à une permutation dont la structure cyclique est:

$$\alpha_1 = \alpha_2 = ... = \alpha_{n-1} = 0; \quad \alpha_n = 1$$

Le nombre de voisinages de n choses correspond à l'ordre de la classe de $S_n$ avec un cycle de longueur n.

La formule de Cauchy donne dans ce cas:

$$\text{ordre de la classe de } S_n \left(\alpha_n = 1\right) = \frac{n!}{n^1 1!}$$

Et ceci n'est rien d'autre que le résultat de Leibniz.

Les partitions: *discerptiones possibiles*

A la fin de sa vie, dans sa correspondance avec Clarke[13], Leibniz expose sa conception de l'espace comme un système de relations et l'arbre généalogique est un "modèle" de cette idée. Presque un demi–siècle auparavant, le jeune penseur se trouvait aux prises avec un aspect combinatoire de la généalogie représentée de manière identique: *Ex hoc ipso problemata origo est numeri personarum in singulis gradibus Arboris Consanguinitatis*[14] ... Il faut donc calculer le "nombre" de *personae*. Ce nombre est donné par $2^n \cdot N$ où n est le "gradus" et N = n + 1 les *cognationes* dans l'arbre. Ici toutefois ce n'est pas ce calcul qui nous retiendra, mais plutôt le versant spatial du problème: assigner une place à une personne donnée. *Persona eo loco intelligantur, ubi puncta sunt*. Ces points de l'espace sont définis par un système de coordonnées: deux nombres –*antecedens* et *sequens*– sont nécessaires pour déterminer ces *puncta* de l'arbre généalogique. Dans ce système, l'ordre des coordonnées importe car (a,b) et (b,a) correspondent à deux points différents. Il s'ensuit une discussion à propos du nombre N de *cognationes*: quel est le nombre des partitions (*discerptiones*) de N en deux parties?

Nous rappellerons qu'une partition de N est une séquence finie non–croissante d'entiers positifs dont la somme est égale à N:

$N_1 + N_2 + ... + N_r = N$; les $N_i$ sont les parties de la partition. D'autre part p(N) est le nombre des partitions de N (par définition p(0) = 1).

Par exemple les partitions de 6 sont:
6; 5,1; 4,2; 4,1,1; 3,3; 3,2,1; 3,1,1,1; 2,2,2; 2,2,1,1; 2,1,1,1,1; 1,1,1,1,1,1.
Et p(6) = 11.

Leibniz établit que le nombre des partitions de N en deux parties est $\frac{N}{2}$ si N est pair et $\frac{N-1}{2}$ si N est impair[15].

---

[13] G. W. Leibniz, "Die philosophische Schriften von G. W. Leibniz", Hrsg. V. C. I. Gerhardt, Berlin, Halle, 1849-1863, vol. 7, p.401.
[14] G. W. Leibniz, A, VI, 1, pp.207-201.
[15] Pour le contexte historique du problème et les contributions leibniziennes, voir Knobloch (On Combinatorics 1974 section 3). D'autre part, dans les articles Déterminants 2001 et Zur Eliminationstheorie 1974 un rapport remarquable entre la théorie des Déterminants et l'arbre généalogique est mis au jour: la virtuosité formelle de Leibniz s'y déploie en toute liberté.

Il nous intéresse de souligner que le nombre des classes de conjugaison du groupe symétrique $S_n$ est égal au nombre des partitions de n.

D'une façon voilée –puisque Leibniz ignorait cette propriété– on décèle la présence de $S_n$ dans la *Dissertatio*.

Une digression s'impose à propos des partitions. Leibniz paraît mesurer toute la difficulté présentée par le calcul de p(N): *sed qui additione datum numerum producendi varietates ... mihi notus non est* et un peu plus loin il parle de "l'abîme des partitions": *At ubi plures partes admittuntur, ingens panditur abyssus discerptionun ....* La signification de cette expression devient plus claire si l'on tient compte du fait que p(N) augmente fort rapidement: p(10) = 22, p(20) = 627, p(100) = 190569292. Dans le premier tiers du XX$^{\text{e}}$ siècle fut trouvée une formule non récurrente pour p(N). Le philosophe d'ailleurs clôt la discussion avec une phrase pleine de bon sens: ce n'est ni l'endroit ni le moment d'en débattre en profondeur (*Exsequi vero hujus loci fortasse, tempori autem non est*).

## IV

Douze problèmes constituent le corpus de la *Dissertatio* dont la structure peut être scindée en trois parties: (a) les problèmes portant sur les complexions 1–3; (b) les problèmes visant les variations de l'ordre et de la disposition: 4–6; (c) les problèmes issus du concept de facteur invariant: 7–12.

Le facteur invariant (*caput variationis*) est défini comme suit: "le facteur invariant est la position des parties certaines; la forme d'une variation de toutes, est celle qui s'obtient de beaucoup de variations, voir plus bas prob. 7".[16] (*Caput variationis est positio certarum partium; Forma variationis, omnium, quæ in pluribus variationibus obtinet, v. infr. prob. 7.*)

Le facteur invariant est une partie d'une variation, présente en toutes les variations. Il s'agit donc d'un facteur avec lequel on doit compter en toute occasion, permanent si l'on peut dire (d'où la pertinence de la traduction).

Leibniz distingue plusieurs facteurs invariants: "en effet les cas sont divers. Car le facteur invariant des variations contient ou bien une chose ou bien beaucoup: si c'est une il est Monadique; ou bien parmi les choses qui doivent varier l'on trouve une autre ou plusieurs autres du même genre que celle contenue dans le facteur invariant. Mais, au contraire le facteur invariant peut comprendre beaucoup de choses, soit que dans le facteur invariant considéré l'on trouve des choses du même genre ou non, à tour de rôle; de même, certaines choses extérieures qui puissent être du même genre que les intérieures ou non".[17]

---

[16] G. W. Leibniz, A, VI, 1, p.173.

Avant d'analyser les cas envisagés, il convient de préciser la position du facteur invariant dans la séquence des thèmes traités. Les complexions –c'est-à-dire l'aspect de la Combinatoire lié à l'arithmétique pure– sont le sujet des problèmes 1–3. Les problèmes 4–6 appartiennent au domaine de l'Arithmétique figurée car ils considèrent des variations d'ordre. Parvenus à ce point, il est bon de rappeler la phrase qui ouvre la *Dissertatio*: "L'Arithmétique est le lieu de cette doctrine. Or, les complexions appartiennent à l'Arithmétique pure, le lieu appartient à l'Arithmétique figurée"[18].

En conséquence, après le sixième problème, le programme annoncé est accompli et, partant, le texte semble conclure C'est alors que l'auteur introduit de manière inattendue la notion de *caput variationis*. Quel sens et quelle portée doit-on attribuer à ce tournant conceptuel de l'œuvre?

Leibniz lui-même répond à ces questions à la fin du Septième problème: "Ce problème put être mené à terme après maintes tentatives fort laborieuses, et sa solution nous demanda beaucoup de travail et de temps. Toutefois personne ne pourrait résoudre d'une autre manière la séquence des problèmes depuis les principes de l'art. En ceux–là, donc, apparaîtra l'usage de celui–ci"[19].

Ce problème dont l'énoncé est concis à l'extrême –Étant donné le facteur invariant trouver les variations– permet de résoudre tous les problèmes précédents, ayant trait aux complexions et aux variations de position, et ce à partir des principes de l'art. Ainsi, à la dérobée, ou presque, une proposition cruciale est avancée: la notion de *caput variationis* est la possibilité même des problèmes formulés et s'érige, par là, en principe de l'Art Combinatoire.

Ainsi en parcourant le problème dix, on lit: "Toute complexion ou variation proposée, moindre et de ces choses, c'est–à–dire celle qui est toute entière contenue dans une autre, est facteur invariant". *(omnis complexio aut variatio propositâ minor et earundem rerum, seu quæ tota in altera continetur, est caput)*

Par conséquent les problèmes 1–6 sont envisagés et résolus en fonction du problème sept qui les rend possibles et leur donne un sens.

Parmi les différents cas l'on distingue des complexions et des variations de disposition. Si le facteur invariant est monadique et fixe, il s'agit d'une variation de l'ordre.

---

[17] G. W. Leibniz, A, VI, 1, pp.219-220..
[18] G.W. Leibniz, A, VI, 1, p.168.
[19] G.W. Leibniz, A, VI, 1, p.221.

Tel est le cas traité en premier lieu: "Le facteur invariant de la variation restant fixe on énumère les choses extérieures et on cherche la variation de celles–ci entre elles (au cas même où elles seraient discontiguës, c'est-à-dire au cas où le facteur invariant se trouverait parmi les variations) se multipliant par ledit facteur, d'après le problème 4, le produit s'appellera A".

Ici, le facteur invariant devient un élément de l'ensemble et le procédé du problème quatre est toujours valable (Pour mémoire, notons que lors de ce problème, Leibniz s'attaque à la question: le nombre des choses étant donné, trouver les variations de l'ordre. Il obtient le résultat $n!$ où $n$ est le nombre de choses qui varient et $n! = n \cdot (n-1) \cdot (n-2) \cdots 3 \cdot 2 \cdot 1$).

Pour ce qui est des complexions, la voie est moins dégagée car par l'entremise du facteur invariant Leibniz ne parvient à rendre compte que du cas des complexions d'un exposant déterminé. Les complexions *simpliciter* échappent à l'emprise de cette notion–clef.

Les complexions de chaque exposant sont associées aux possibles facteurs invariants homogènes avec les choses qui peuvent varier.

Par exemple si les choses sont {a, b, c, d, e} et le facteur invariant, homogène avec tous les membres extérieurs, est ab, alors le calcul de toutes les com2naisons possibles (voir Glossaire) dans cet ensemble détermine tous les facteurs invariants possibles.

Cette méthode peut s'appliquer aux facteurs invariants, quel que soit l'exposant fixé. Par suite, les complexions *simpliciter* sont les seules qui demeurent hors d'atteinte. Aussi peut–on signaler que les complexions et les variations d'ordre se réduisent à des modes engendrés par le facteur invariant: la variation d'ordre surgit lorsque le facteur invariant remplit certaines conditions et à leur tour, les complexions correspondent aux déterminations dudit facteur.

Si l'on s'en tient à ce rapport, on peut conclure que le domaine de l'Arithmétique pure et le domaine de l'Arithmétique figurée deviennent, en s'assemblant, une troisième contrée dont il faut préciser la nature.

L'analyse antérieure a mis à jour la racine conceptuelle de l'œuvre: si tout d'abord la *Dissertatio* se présentait sous la forme d'un traité sur les complexions et les variations, maintenant elle se révèle comme une théorie du facteur invariant.

Il est à présent loisible d'apprécier la portée de deux énoncés du philosophe, effleurés plus haut:

- "Le facteur invariant est la position des parties certaines" (*positio certarum partium*).

- "On dit d'une chose qu'elle est homogène s'il est possible de la mettre dans une place donnée, sauf dans le facteur invariant"[20].

Dans ces deux cas un trait du concept saute aux yeux: le facteur invariant est une "caractéristique spatiale" (ci–après, par abus de langage, nous dirons: un espace).

En outre il s'agit d'un espace bien défini qui organise tous les autres espaces. L'expression latine –*caput variationis*– décrit avec justesse la fonction de cet élément qui confère un ordre et un sens au Tout. La *Dissertatio* fonde le calcul de tous les occupants possibles de cet espace (complexions) et de tous les ordres possibles du monde issus du facteur invariant. L'Art Combinatoire rompt les cadres de l'algorithme arithmétique et va au-delà de la Combinatoire. En fait son point de fuite et son point d'appui coïncident: "parce que l'on verra sortir toutes choses du plus profond de la doctrine des variations, laquelle, toute seule de fait, conduit l'âme docile à travers le tout infini et comprend tout ensemble l'harmonie du monde et les constructions ultimes des choses et la série des formes, dont l'incroyable utilité sera estimée, à son juste prix, –en dernier lieu– par la parfaite ou presque parfaite philosophie."[21]

La distinction Arithmétique pure – Arithmétique figurée vise à établir deux étapes dans le chemin qui mène vers une plus haute discipline. L'idée de *figure* élucide cette différence.

Il y a lieu de penser que Leibniz reprend la conception de la quantité représentée dans l'espace, attribuée à l'école de Pythagore. Le philosophe de Leipzig rattache l'arithmétique figurée aux variations d'ordre car celles–ci portent sur la totalité des éléments et n'envisagent que les changements de leur distribution spatiale. Par contre les complexions considèrent une partie du tout sans tenir compte de la variation de lieu.

Résumons les étapes de notre exposé:
i)   La *Dissertatio* est la doctrine du facteur invariant.
ii)  Le facteur invariant est un lieu (une place, un endroit).
iii) Les complexions déterminent les occupants possibles dudit lieu.
iv)  Les variations d'ordre définissent toutes les dispositions possibles des objets à partir de tous les facteurs invariants possibles.

La proposition suivante condense ces résultats:

*La doctrine du facteur invariant s'incarne dans le calcul et la détermination de toutes les figurations possibles de l'Univers.*

---

[20] G.W.Leibniz, A, VI, 1, p.173.
[21] G.W.Leibniz, A, VI, 1, p.187.

Ici *figuration* signifie *projection du nombre dans l'espace*.

Pour Leibniz la quantité –détermination d'une totalité abstraite– est immatérielle. Le nombre n'est point mesure mais unité. Donc le nombre concerne tous les étants, non seulement ceux qui sont corporels. Si l'on borne l'attribution du nombre à un certain type de réalité on restreint d'autant le champ de la découverte puisque le nombre parachève une totalité nouvelle. L'Univers naît de la composition d'unités. Rien n'est au-delà du nombre. Tout est nombre.

A oublier ces considérations on perd les nuances et la portée du verbe *calculer* chez Leibniz dont la pensée est faite de finesse conceptuelle poussée à l'extrême. Aussi ne choisit–il jamais l'une des figurations: n'importe quel étant à n'importe quel moment peut occuper la place du facteur invariant et devenir le centre du monde.

Dans ce contexte, les complexions et les variations sont les instruments d'un algorithme majeur: le Calcul de toutes les possibilités de composer l'Univers. Cet enjeu pointe à l'horizon de chacune des pages de la *Dissertatio*. Le facteur invariant est un lieu. Ipso facto tous les autres facteurs, à leur tour, deviennent des lieux. Les différents ordres du monde sont ainsi conçus comme des ordres spatiaux. D'où l'importance de l'idée de figure et le primat des variations d'ordre. Nulle hiérarchie parmi les étants: le nombre leur confère unité et à travers lui paraît l'harmonie de l'Univers.

L'Art Combinatoire s'évertue à trouver le nombre du monde projeté en l'espace.

### Glossaire

1.– <u>Complexion</u>: correspond à la notion actuelle de Combinaison, c'est–à–dire à la détermination de tous les sous–ensembles possibles ayant le même nombre d'éléments qu'un ensemble donné. Dans la terminologie de la *Dissertatio* un sous–ensemble est un tout plus petit compris dans un tout plus grand et un élément est appelé partie. Leibniz considère la combinaison comme un cas spécial de Complexion pour lequel le nombre d'éléments de chaque sous–ensemble est deux. Dans sa notation, combinaison = com2naison. Si le nombre des éléments est 3, 4, etc. il écrit, respectivement, com3naison, com4naison, etc. Les complexions d'un seul élément se dénomment Unions.

2.– <u>Ordre</u>: a) équivaut au concept de cardinalité, c'est–à–dire au nombre d'éléments (parties) d'un ensemble (tout). Sur ce point, toutefois, Leibniz est ambigu car il envisage aussi

b) les Variations d'ordre qui ne concernent pas la quantité des parties mais les changements de la disposition d'une quantité donnée d'avance de parties. Ainsi <u>abc</u> varie en devenant <u>bca</u>.

3.– <u>Variation</u>: d'après la définition il s'agit d'un Changement de Relation. Or un tel Changement advient de trois manières: Rapport, Lieu et Conjonction. Et porte sur trois niveaux: substantiel, qualitatif et quantitatif.

Tout ce qui existe est composé et résulte des rapports des parties composantes. Par suite un changement de Relation peut être:

a) un changement en l'essence d'une chose.
b) un changement par rapport au lieu: les parties changent de place sans modifier l'essence de la chose mais en suscitant des changements qualitatifs ou quantitatifs.
c) un changement issu de la relation d'une chose avec une autre chose. Cette opération s'appelle conjonction et il peut s'ensuivre un être nouveau ou bien un changement moins radical, par exemple de la figure, de la quantité, de la couleur etc.

4.– <u>Complexion</u> *Simpliciter*: cette expression dénote la somme de toutes les complexions possibles en un tout. À titre d'illustration soit un tout (ensemble) composé de quatre parties (éléments): les complexions *simpliciter* résultent de la somme des Unions, des com2naisons, des com3naisons et des com4naisons.

5.– <u>Variation de Disposition</u>: ou variation de voisinage. Ce concept porte sur le changement de la relation spatiale entre les parties. C'est par conséquent un changement par rapport à un espace relatif. Si l'on dispose trois parties <u>a</u>, <u>b</u>, <u>c</u> en cercle, <u>abc</u> et <u>bca</u> sont un même état des choses. Par contre le passage de <u>abc</u> à <u>cba</u> détermine un changement de voisinage.